\documentclass[12pt,amstex,amssymb,a4paper]{article}
\usepackage{amstex}
\usepackage{amssymb}
\newtheorem{theorem}{Theorem}%[section]
\newtheorem{lemma}[theorem]{Lemma}
\newtheorem{proposition}[theorem]{Proposition}

\newcommand{\prf}{{\em Proof}. }
\newcommand{\qed}{\hspace*{\fill}$\Box$}

\newcommand{\lra}{\longrightarrow} 
\newcommand{\Buparrow}{\Big\uparrow}
\newcommand{\ka}{{\cal A}}
\newcommand{\kb}{{\cal B}}

\newcommand{\kk}{{\cal K}}
\newcommand{\kl}{{\cal L}}
\newcommand{\kn}{{\cal N}}
\newcommand{\ko}{{\cal O}}
\newcommand{\kt}{{\cal T}}
\newcommand{\IA}{{\mathbb A}}
\newcommand{\IC}{{\mathbb C}}
 
\newcommand{\IP}{{\mathbb P}}

\newcommand{\dual}{\makebox[0mm]{}^{{\scriptstyle\vee}}}
\newcommand{\tensor}{\otimes}

\newcommand{\isom}{\cong}

\newcommand{\beeq}[1]{\begin{eqnarray}\label{#1}}
\newcommand{\eneq}{\end{eqnarray}}

\newcommand{\Hilb}{{\rm Hilb}}
\newcommand{\Quot}{{\rm Quot}}
\newcommand{\Supp}{{\rm Supp}}

\newcommand{\Hom}{{\rm Hom}}
\newcommand{\Ext}{{\rm Ext}}
\newcommand{\Socle}{{\rm Soc}}
\newcommand{\id}{{\rm id}}
\newcommand{\im}{{\rm im}}
\newcommand{\Ker}{{\rm ker}}
\newcommand{\Cok}{{\rm coker}}
\newcommand{\codim}{{\rm codim}}
\newcommand{\ses}[3]{0\rightarrow#1\rightarrow#2\rightarrow#3\rightarrow0}
 
\begin{document}
\title{Irreducibility of the Punctual Quotient Scheme of a Surface}
\author{Geir Ellingsrud\and Manfred Lehn}
\date{January 22, 1997}
\maketitle

Let $X$ be a smooth projective surface, $E$ a locally free sheaf of rank
$r\geq 1$ on $X$, and let $\ell\geq 1$ be an integer.
$\Quot(E,\ell)$ denotes Grothendieck's quotient scheme \cite{Grothendieck} that
parametrises all surjections $E\to T$, where $T$ is a zero-dimensional sheaf of 
length $\ell$, modulo automorphisms of $T$.
Sending a quotient $E\to T$ to the point
$\sum_{x\in X}\ell(T_x)x$ in the symmetric product $S^\ell(X)$ defines
a morphism $\pi:\Quot(E,\ell)\to S^\ell(X)$ \cite{Grothendieck}.
It is the purpose of this note to prove the following theorem:

\begin{theorem}--- $\Quot(E,\ell)$ is an irreducible scheme of dimension
$\ell(r+1)$. The fibre of the morphism $\pi:\Quot(E,\ell)\to S^{\ell}(X)$
over a point $\sum_x \ell_xx$ is irreducible of dimension $\sum_x (r\ell_x-1)$.
\end{theorem}

If $r=1$, i.e.\ if $E$ is a line bundle, then $\Quot(E,\ell)$ is isomorphic
to the Hilbert scheme $\Hilb^\ell(X)$. For this case, the first assertion of
the theorem is due to Fogarty \cite{Fogarty}, whereas the second assertion was
proved by Brian\c{c}on \cite{Briancon}. For general $r\geq 2$, the first 
assertion of the theorem is a result due to J.\ Li  and D.\ Gieseker 
\cite{Li},\cite{GiesekerLi}. We give a different proof with a more geometric
flavour, generalising a technique from Ellingsrud and Str{\o}mme \cite{EllingsrudStromme}. The second assertion is a new result for $r\geq 2$.

\section{Elementary Modifications}\label{modific}

Let $X$ be a smooth projective surface and $x\in X$. 
If $N$ is a coherent $\ko_X$-sheaf, $e(N_x)=\hom_X(N,k(x))$ denotes
the dimension of the fibre $N(x)$, which by Nakayama's Lemma is the
same as the minimal number of generators of the stalk $N_x$. If $T$ is a
coherent sheaf with zero-dimensional support, we denote
by $i(T_x)=\hom_X(k(x),T)$ the dimension of the
socle of $T_x$, i.e.\ the submodule $\Socle(T_x)\subset T_x$ of all elements 
that are annihilated by the maximal ideal in $\ko_{X,x}$.

\begin{lemma}\label{comparisonofiande}--- 
Let $[q:E\to T]\in\Quot(E,\ell)$ be a closed point and let
$N$ be the kernel of $q$. Then the socle dimension of $T$ and the number of
generators of $N$ at $x$ are related as follows:
$$e(N_x)=i(T_x)+r.$$
\end{lemma}

\prf Write $e(N_x)=r+i$ for some integer $i\geq 0$. Then 
there is a minimal free resolution
$0\lra\ko_{X,x}^i\stackrel{\alpha}{\lra}\ko_{X,x}^{r+i}\lra N_x\lra0$,
where all coefficients of the homomorphism
$\alpha$
are contained in the maximal ideal of $\ko_{X,x}$. We have $\Hom(k(x),T_x)\isom\Ext_X^1(k(x),N_x)$ and applying the functor
$\Hom(k(x),\,.\,)$ one finds an
exact sequence
$$0\lra\Ext_X^1(k(x),N_x)\lra\Ext_X^2(k(x),\ko_{X,x}^i)
\stackrel{\alpha'}{\lra}\Ext_X^2(k(x),\ko_{X,x}^{r+i}).$$
But as $\alpha$ has coefficients in
the maximal ideal, the homomorphism $\alpha'$ is zero. Thus $\Hom(k(x),T)\isom\Ext_X^2(k(x),\ko_{X,x}^i)\isom k(x)^i$.
\qed

The main technique for proving the theorem will be induction on the length
of $T$. Let $N$ be the kernel of a surjection $E\to T$, let $x\in X$ be
a closed point, and let $\lambda:N\to k(x)$ be any surjection.
Define a quotient $E\to T'$ by means of the following push-out diagram:
$$\begin{array}{ccccccccc}
&&0&&0\\
&&\uparrow&&\uparrow\\
0&\lra&k(x)&\stackrel{\mu}{\lra}&T'&\lra&T&\lra&0\\
&&{\scriptstyle\lambda}\uparrow{\phantom{\scriptstyle\lambda}}&&\uparrow&&
\|\\
0&\lra&N&\lra&E&\lra&T&\lra&0\\
&&\uparrow&&\uparrow&&\\
&&N'&=&N'\\
&&\uparrow&&\uparrow\\
&&0&&0 
\end{array}
$$
In this way every element $\langle\lambda\rangle\in\IP(N(x))$ determines a 
quotient $E\to T'$ together with an element $\langle\mu\rangle\in
\IP(\Socle(T'_x)\dual)$. (Here $W\dual:=\Hom_k(W,k)$ denotes the vector space  dual to $W$.) Conversely, if $E\to T'$ is given, any such 
$\langle\mu\rangle$ determines $E\to T$ and a point $\langle\lambda\rangle$.
We will refer to this situation by saying that $T'$ is obtained from $T$
by an elementary modification.

We need to compare the invariants for $T$ and $T'$: Obviously, $\ell(T')=
\ell(T)+1$. Applying the functor $\Hom(k(x),\,.\,)$ to the upper row in the
diagram we get an exact sequence
$$0\longrightarrow k(x)\longrightarrow\Socle(T'_x)\to\Socle(T_x)\longrightarrow
\Ext_X^1(k(x),k(x))\isom k(x)^2,$$
and therefore $|i(T_x)-i(T'_x)|\leq 1$. Moreover, 
we have $0\leq e(T'_x)-e(T_x)\leq 1$. Two cases deserve
more attention:

\begin{lemma}\label{caseeincreases}--- Consider the natural 
homomorphisms $g:N(x)\to E(x)$ and $f:\Socle(T'_x)\to T'\to T'(x)$.
The following assertions are equivalent
\begin{enumerate}
\item $e(T')=e(T)+1$
\item $\langle\mu\rangle\not\in\IP(\Ker(f)\dual)$
\item $\langle\lambda\rangle\in\IP(\im (g))$.
\end{enumerate}
Moreover, if these conditions are satisfied, then 
$T'\isom T\oplus k(x)$ and $i(T'_x)=i(T_x)+1$.
\end{lemma}

\prf Clearly, $e(T')=e(T)+1$ if and only if $\mu(1)$ represents a non-trivial
element in $T'(x)$ if and only if $\mu$ has a left inverse if and only if
$\lambda$ factors through $E$.\qed

\begin{lemma}\label{caseiincreases}--- Still keeping the notations above, let $E\to T'_{\lambda}$ be the modification of $E\to T$ determined by the point
$\langle\lambda\rangle\in\IP(N(x))$. Similarly, for $\langle\mu'\rangle\in
\IP(\Socle(T_x)\dual)$ let $T^-_{\mu'}=T/\mu'(k(x))$.
If $i(T'_{\lambda,x})=i(T_x)+1$ for all $\langle\lambda\rangle\in\IP(N(x))$,
then $i(T_x)=i(T^-_{\mu',x})-1$ for all $\langle\mu'\rangle\in\IP(\Socle(T_x)\dual)$ as well.
\end{lemma}

\prf Let $\Phi:\Hom_X(N,k(x))\to\Hom_k(\Ext_X^1(k(x),N),\Ext_X^1
(k(x),k(x)))$ be the homomorphism which is adjoint to the natural pairing
$$\Hom_X(N,k(x))\tensor\Ext_X^1(k(x),N)\to \Ext_X^1(k(x),k(x)).$$
Identifying $\Socle(T_x)\isom\Ext_X^1(k(x),N)$, we see that $i(T'_{\lambda,x})=1+i(T_x)-{\rm rank}(\Phi(\lambda))$. 
The action of $\Phi(\lambda)$ on a socle element $\mu':k(x)\to T$ can be described by the following diagram of pull-backs and push-forwards
$$\begin{array}{ccccccccc}
0&\to&N&\to&E&\to&T&\to&0\\
&&\|&&\uparrow&&\phantom{{\scriptstyle\mu}}\uparrow{\scriptstyle\mu'}\\
0&\to&N&\to&N^-_{\mu'}&\to&k(x)&\to&0\\
&&{\scriptstyle\lambda}
\downarrow\phantom{{\scriptstyle\lambda}}&&\downarrow&&\|\\
0&\to&k(x)&\to&\xi&\to&k(x)&\to&0
\end{array}$$
The assumption that $i(T'_{\lambda,x})=1+i(T_x)$ for all $\lambda$, is equivalent to $\Phi=0$.  This implies that for every $\mu'$ and every $\lambda$ the extension in the third row splits, which in turn means that every $\lambda$
factors through $N^-_{\mu'}$, i.e.\ that $N(x)$ embeds into $N^-_{\mu'}(x)$. 
Hence, for $T^-_{\mu'}=E/N^-_{\mu'}=\Cok(\mu)$ we get
$i(T^-_{\mu',x})=e(N^-_{\mu',x})-r=e(N_x)+1-r=i(T_x)+1$.\qed

\section{The Global Case}\label{globalcasesection}

Let $Y_\ell=\Quot(E,\ell)\times X$, and consider the universal exact sequence
of sheaves on $Y_\ell$:
$$\ses{\kn}{\ko_{\Quot}\tensor E}{\kt}.$$
The function $y=(s,x)\mapsto i(\kt_{s,x})$ is upper semi-continuous. Let
$Y_{\ell,i}$ denote
the locally closed subset $\{y=(s,x)\in Y_\ell|i(\kt_{s,x})=i\}$ with the 
reduced subscheme structure.

\begin{proposition}\label{globprop}--- $Y_\ell$ is irreducible of
dimension $(r+1)\ell+2$.
For each $i\geq 0$ one has $\codim(Y_{\ell,i},Y_\ell)\geq 2i$,
\end{proposition}

Clearly, the first assertion of the theorem follows from this.

\prf The proposition will be proved by induction on $\ell$, the case $\ell=1$
being trivial: $Y_1=\IP(E)\times X$, the stratum $Y_{1,1}$ is the graph of the 
projection $\IP(E)\to X$ and $Y_{1,i}=\emptyset$ for $i\geq 2$. 
Hence suppose the proposition has been proved for some $\ell\geq 1$.

We describe the `global' version of the elementary modification discussed
above. Let $Z=\IP(\kn)$ be the projectivization of the family $\kn$ and let $\varphi=(\varphi_1,\varphi_2):Z\to Y_\ell=\Quot(E,\ell)\times X$ denote the
natural projection morphism. On $Z\times X$ there is canonical epimorphism
$$\Lambda:(\varphi_1\times\id_X)^*\kn\to(\id_Z,\varphi_2)_*\varphi^*\kn\to
(\id_Z,\varphi_2)_*\ko_{Z}(1)=:\kk.$$
As before we define a family $\kt'$ of quotients of length $\ell+1$ by means of $\Lambda$:
$$\begin{array}{ccccccccc}
0&\lra&\kk&\lra&\kt'&\lra&(\varphi_1,\id_X)^*\kt&\lra&0\\
&&\Lambda\Buparrow\phantom{\Lambda}&&\Buparrow&&\Big\|\\
0&\lra&(\varphi_1,\id_X)^*\kn&\lra&\ko_Z\tensor E&\lra&(\varphi_1,\id_X)^*\kt&
\lra&0
\end{array}$$
Let $\psi_1:Z\to \Quot(E,\ell+1)$ be the classifying morphism for the family
$\kt'$, and define
$\psi:=(\psi_1,\psi_2:=\varphi_2):Z\to Y_{\ell+1}$.
The scheme $Z$ together with the morphisms $\varphi:Z\to Y_\ell$ and $\psi:
Z\to Y_{\ell+1}$ allows us to relate the strata $Y_{\ell,i}$ and $Y_{\ell+1,j}$.
Note that $\psi(Z)=\bigcup_{j\geq1}Y_{\ell+1,j}$.

The fibre of $\varphi$ over a point $(s,x)\in Y_{\ell,i}$ is given by $\IP(\kn_s(x))\isom\IP^{r-1+i}$, since $\dim(\kn_s(x))=r+i(T_{s,x})=r+i$ by
Lemma \ref{comparisonofiande}.
Similarly, the fibre of $\psi$ over a point $(s',x)\in Y_{\ell+1,j}$ is 
given by $\IP(\Socle(\kt'_{s',x})\dual)\isom\IP^{j-1}$.
If $T'$ is obtained from $T$ by an elementary modification, then $|i(T')-i(T)|
\leq 1$ as shown above. This can be stated in terms of $\varphi$ and $\psi$
as follows: For each $j\geq 1$ one has:
$$\psi^{-1}(Y_{\ell+1,j})\subset\bigcup_{|i-j|\leq1}\varphi^{-1}(Y_{\ell,i}).$$
Using the induction hypothesis on the dimension of $Y_{\ell,i}$ and the 
computation of the fibre dimension of $\varphi$ and $\psi$, we get
$$\dim(Y_{\ell+1,j})+(j-1)\leq\max_{|i-j|\leq 1}\{(r+1)\ell+2-2i+(r-1+i)\}$$
and
$$\dim(Y_{\ell+1,j})\leq(r+1)(\ell+1)+2-2j-\min_{|i-j|\leq 1}\{i-j+1\}.$$
As $\min_{|i-j|\leq 1}\{i-j+1\}\geq0$, this proves the dimension estimates
of the proposition.

It suffices to show that $Z$ is irreducible.
Then $\Quot(E,\ell+1)=\psi_1(Z)$ and $Y_{\ell+1}$ are irreducible as well. 

Since $X$ is a smooth surface, the epimorphism
$\ko_{\Quot}\tensor E\to \kt$ can be completed to a finite resolution
$$0\longrightarrow\ka\longrightarrow\kb\longrightarrow\ko_{\Quot}\tensor E
\longrightarrow \kt\longrightarrow0$$
with locally free sheaves $\ka$ and $\kb$ on $Y_{\ell}$ of rank $n$ and $n+r$,
respectively, for some positive integer $n$.
It follows that $Z=\IP(\kn)\subset \IP(\kb)$ is the vanishing locus of the
composite homomorphism $\varphi^*\ka\to\varphi^*B\to\ko_{\IP(\kb)}(1)$.
In particular, assuming by induction that $Y_\ell$ is irreducible, $Z$ is
locally cut out from an irreducible variety of dimension $(r+1)\ell+2+(r+n-1)$
by $n$ equations. Hence every irreducible component of $Z$ has dimension at
least $(r+1)(\ell+1)$. But the dimension estimates for the stratum $Y_{\ell,i}$
and the fibres of $\varphi$ over it yield:
$$\dim(\varphi^{-1}(Y_{\ell,i}))\leq (r+1)\ell+2-2i+(r+i-1)=(r+1)(\ell+1)-i,$$
which is strictly less than the dimension of any possible component of $Z$, if
$i\geq 1$. This implies that the irreducible variety $\varphi^{-1}(Y_{\ell,0})$
is dense in $Z$. Moreover, since the fibre of $\psi$ over $Y_{\ell+1,1}$
is zero-dimensional, $\dim(Y_{\ell+1})=\dim(Y_{\ell+1,1})+2=\dim(Z)+2$
has the predicted value.\qed

\section{The Local Case}\label{localcasesection}

We now concentrate on quotients $E\to T$, where $T$ has support in a single
fixed closed point $x\in X$. For those quotients the structure of $E$ is of
no importance, and we may assume that $E\isom\ko^r_X$. Let $Q^r_\ell$ denote
the closed subset
$$\Big\{[\ko^r_X\to T]\in\Quot(\ko^r_X,\ell)|\,\,\Supp(T)=\{x\}\Big\}$$ with
the reduced subscheme structure. We may consider
$Q^r_\ell$ as a subscheme of ${Y_{\ell,1}}$ by sending $[q]$ to $([q],x)$.
Then it is easy to see that $\varphi^{-1}(Q^r_\ell)=\psi^{-1}(Q^r_{\ell+1})$. 
Let this scheme be denoted by $Z'$. 

We will use a stratification of $Q^r_\ell$ both by the socle dimension $i$ and
the number of generators $e$ of $T$ and denote the corresponding locally closed
subset by $Q^{r,e}_{\ell,i}$. Moreover, let
$Q^{r}_{\ell,i}=\bigcup_{e}Q^{r,e}_{\ell,i}$ and define $Q^{r,e}_\ell$ similarly.
Of course, $Q^{r,e}_{\ell,i}$ is empty unless $1\leq i\leq \ell$ and $1\leq e \leq \min\{r,\ell\}$.

To prove the second half of the theorem it suffices to show:

\begin{proposition}\label{localcase}--- $Q^r_\ell$ is an irreducible variety of dimension $r\ell-1$.
\end{proposition}

\begin{lemma}\label{dimensionestimatesforstrata}--- $\dim(Q^{r,e}_{\ell,i})
\leq(r\ell-1)-(2(i-1)+\binom{e}{2}).$
\end{lemma}

\prf By induction on $\ell$: if $\ell=1$, then $Q^r_1\isom\IP^{r-1}$, and
$Q^{r,e}_{1,i}=\emptyset$ if $e\geq 2$ or $i\geq 2$. Assume that the lemma has 
been proved for some $\ell\geq 1$.

Let $[q':\ko_X^r\to T']\in Q^{r,e}_{\ell+1,j}$ be a closed point. Suppose
that the map $\mu:k(x)\to T'(x)$ represents a point in
$\psi^{-1}([q'])=\IP(\Socle(T'_x)\dual)$ and that $T_\mu=\Cok(\mu)$ is the corresponding modification. If $i=i(T_{\mu,x})$ and $\varepsilon=e(T_{\mu,x})$, then, according to Section \ref{modific}, the pair $(i,\varepsilon)$ can take the following values:
\beeq{fourpossibilities}
(i,\varepsilon)=(j-1,e-1),\,(j-1,e),\,(j,e)\,\mbox{ or }\, (j+1,e),
\eneq
in other words:
$$\psi^{-1}(Q^{r,e}_{\ell+1,j})\subset
\varphi^{-1}(Q^{r,e-1}_{\ell,j-1})\cup\bigcup_{|i-j|\leq1}
\varphi^{-1}(Q^{r,e}_{\ell,i}).$$
Subdivide $A=Q^{r,e}_{\ell,j}$ into four locally closed subsets $A_{i,\varepsilon}$ according to the generic value of $(i,\varepsilon)$ on
the fibres of $\psi$. Then
$$\dim(A_{i,\varepsilon})+(j-1)\leq \dim(Q^{r,\varepsilon}_{\ell,i})+d_{i,\varepsilon},$$
where $d_{i,\varepsilon}$ is the fibre dimension of the morphism
$$\varphi:\psi^{-1}(A_{i,\varepsilon})\cap\varphi^{-1}(Q^{r,\varepsilon}_{\ell,i})
\lra Q^{r,\varepsilon}_{\ell,i}.$$
By the induction hypothesis we have bounds for $\dim(Q^{r,\varepsilon}_{\ell,i})$, and we can bound $d_{i,\varepsilon}$ in the 
four cases (\ref{fourpossibilities}) as follows:

A) If $[q:\ko_X^r\to T]\in Q^{r,e-1}_{\ell,j-1}$ is a closed point with 
$N=\Ker(q)$, then according to Lemma \ref{caseeincreases}
\begin{eqnarray*}
\varphi^{-1}([q])\cap\psi^{-1}(A_{e-1,j-1})&\isom&\IP(\im(g:N(x)\to k(x)^r))\\
&\isom&\IP(\Ker(k(x)^r\to T(x))\isom \IP^{r-e},
\end{eqnarray*}
since $\im(k(x)^r\to T(x))\isom k^{e-1}$. Hence $d_{j-1,e-1}=r-e$ and 
\begin{eqnarray*}
\dim(A_{j-1,e-1})&\leq&\dim Q^{r,e-1}_{\ell,j-1}+(r-e)-(j-1)\\
&\leq&\left\{(r\ell-1)-2(j-2)-\binom{e-1}{2}\right\}+(r-e)-(j-1)\\
&=&\left\{(r(\ell+1)-1)-2(j-1)-\binom{e}{2}\right\}-(j-2).
\end{eqnarray*}
Note that this case only occurs for $j\geq2$, so that $(j-2)$ is always
nonnegative.

B) In the three remaining cases
$$\varepsilon=e\mbox{ and } i=j-1,\, j, \mbox{ or } j+1$$
we begin with the rough estimate $d_{i,e}\leq r+i-1$ as in Section \ref{globalcasesection}. This yields:
\begin{eqnarray}
\label{rechnung2}
\dim(A_{i,e})&\leq&\left\{(r\ell-1)-2(i-1)-\binom{e}{2}\right\}+(r+i-1)-(j-1)\\
\label{rechnung3}
&=&\left\{(r(\ell+1)-1)-2(j-1)-\binom{e}{2}\right\}-(i-j).
\end{eqnarray}
Thus, if $i=j$ we get exactly the estimate asserted in the Lemma, if $i=j+1$
the estimate is better than what we need by 1, but if $i=j-1$, the estimate
is not good enough and fails by 1. It is this latter case that we must further
study: let $[q:\ko^r_X\to T]$ be a point in $Q^{r,e}_{\ell,j-1}$ with $N=\Ker(q)$. By Lemma
\ref{caseiincreases} there are two possibilities:
\begin{itemize}
\item[---]{\em Either} the fibre
$\varphi^{-1}([q])\cap\psi^{-1}(A_{j-1,e})$ is a {\em proper} 
closed subset of $\IP(N(x))$ which improves the estimate for the dimension of 
the fibre $\varphi^{-1}([q])$ by 1,
\item[---]{\em or} this fibre {\em equals} with $\IP(N(x))$, in which case  we have $i(T^-)=i(T)+1$ for every modification $T^-=\Cok(\mu^-:k(x)\to T)$. 
But, as we just saw, calculation (\ref{rechnung3}), applied to the contribution of $Q^{r,e}_{\ell-1,j}$ to $Q^{r,e}_{\ell,j-1}$, shows that the dimension 
estimate 
for the locus of such points $[q]$ in $Q^{r,e}_{\ell,j-1}$ can be improved by 1 compared to the dimension estimate for $Q^{r,e}_{\ell,j-1}$ as stated in the lemma.
\end{itemize}
Hence in either case we can improve estimate (\ref{rechnung3}) by 1 and get
$$\dim(A_{j-1,e})\leq(r(\ell+1)-1)-2(j-1)-\binom{e}{2}$$
as required.
Thus, the lemma holds for $\ell+1$.\qed

\begin{lemma}\label{intheclosureofelessthanr}--- 
$\psi(\varphi^{-1}(Q^{r,e}_{\ell}))\subset\overline{Q^{r,e}_{\ell+1}}$.
\end{lemma}

\prf Let $[q:\ko^r_X\to T]\in Q^{r,e}_{\ell,i}$ be a closed point with 
$N=\Ker(q)$. Then $\varphi^{-1}([q])=\IP(N(x))\isom\IP^{r+i-1}$ and
$\varphi^{-1}([q])\cap\psi^{-1}(Q^{r,e+1}_{\ell+1})\isom \IP(\im (G))\isom
\IP^{r-e-1}$. Since we always have $e\geq 1,i\geq 1$, a dense open part of
$\varphi^{-1}([q])$ is mapped to $Q^{r,e}_{\ell+1}$.\qed

\begin{lemma}\label{inductionstepforr}--- 
If $r\geq 2$ and if $Q^{r-1}_\ell$ is irreducible of dimension
$(r-1)\ell-1$, then $Q^{r,<r}_{\ell}:=\bigcup_{e<r}Q^{r,e}_{\ell}$ is an
irreducible open subset of $Q^{r}_\ell$ of dimension $r\ell-1$.
\end{lemma}

\prf Let $M$ be the variety of all $r\times(r-1)$ matrices over $k$ of maximal
rank, and let $\ses{\ko_M^{r-1}}{\ko_M^r}{\kl}$ be the corresponding 
tautological sequence of locally free sheaves on $M$. 
Consider the open subset $U\subset M\times Q^{r}_\ell$ of points 
$(A,[\ko^{r}\to T])$ such that the composite homomorphism
$$\ko^{r-1}\stackrel{A}{\lra}\ko^{r}\lra T$$
is surjective. Clearly, the image of $U$ under the projection to
$Q^{r}_\ell$ is $Q^{r,<r}_\ell$.
On the other hand, the tautological epimorphism
$$\ko_{U\times X}^{r-1}\to\ko_{U\times X}^{r}\to(\ko_M\tensor\kt)|_{U\times X}$$ 
induces a classifying morphism $g':U\to Q^{r-1}_\ell$. The morphism
$g=(pr_1,g'):U\to M\times Q^{r-1}_\ell$ is surjective. In fact, it is an affine
fibre bundle with fibre
$$g^{-1}(g(A,[\ko^{r-1}\to T]))\isom\Hom_k(\kl(A),T)\isom \IA_k^\ell.$$
Since $Q^{r-1}_\ell$ is irreducible of dimension $(r-1)\ell-1$ by 
assumption, $U$ is irreducible of dimension $r\ell-1+\dim(M)$,
and $Q^{r,<r}_\ell$ is irreducible of dimension $r\ell-1$.\qed

{\em Proof of Proposition \ref{localcase}.} The irreducibility of 
$Q^{r}_{\ell}$ will be proved by induction over $r$ and $\ell$:
the case $(\ell=1,r
\mbox{ arbitrary })$ is trivial; whereas $(\ell\mbox{ arbitrary },r=1)$ is the
case of the Hilbert scheme, for which there exist several proofs (\cite{Briancon}, \cite{EllingsrudStromme}). 
Assume therefore that $r\geq 2$ and that the proposition holds for
$(\ell,r)$ and $(\ell+1,r-1)$. We will show that it holds for $(\ell+1,r)$ as
well.

Recall that $Z':=\varphi^{-1}(Q^{r}_\ell)=Q^r_{\ell}\times_{Y_{\ell}}Z$.
Every irreducible component of $Z'$ has dimension greater than or
equal to $\dim(Q^{r}_\ell)+r-1=r(\ell+1)-2$ (compare Section \cite{globalcasesection}). On the other hand,
$\dim(\varphi^{-1}(Q^r_{\ell,i}))\leq r\ell-1-2(i-1)+(r+i-1)=r(\ell+1)-i$.
Thus an irreducible components of $Z'$ is either the closure of
$\varphi^{-1}(Q^r_{\ell,1})$ (of dimension $r(\ell+1)-1)$) or the closure of
$\varphi^{-1}(W)$ for an irreducible component $W\subset Q^r_{\ell,2}$ of
maximal possible dimension $r\ell-3$. But according to Lemma 
\ref{intheclosureofelessthanr} the image of $\varphi^{-1}(W)$ under $\psi$
will be contained in the closure of $Q^{r,<r}_{\ell+1}$, unless $W$ is contained
in $Q^{r,r}_{\ell,2}$. But Lemma \ref{dimensionestimatesforstrata} says that
$Q^{r,r}_{\ell,2}$ has codimension $\geq 2+\binom{r}{2}\geq 3$ if $r\geq 2$,
and hence cannot contain $W$ for dimension reasons. Hence any irreducible
component of $Z'$ is mapped by $\psi$ into the closure of $Q^{r,<r}_{\ell+1}$
which is irreducible by Lemma \ref{inductionstepforr} and the induction
hypothesis. This finishes the proof of the proposition.\qed

\end{document}